\documentclass{amsart} 

\usepackage{mathtools}
\usepackage{amsmath}
\usepackage{amsthm}
\usepackage{amssymb}
\usepackage{amsfonts}
\usepackage{enumerate}
\usepackage{algorithm}
\usepackage[noend]{algpseudocode}
\usepackage{appendix}
\usepackage{xcolor}

\newtheorem{theorem}{Theorem}[section]
\newtheorem{lemma}[theorem]{Lemma}

\newtheorem{corollary}[theorem]{Corollary}

\theoremstyle{definition}
\newtheorem{definition}[theorem]{Definition}

\theoremstyle{remark}

\numberwithin{equation}{section}

\newcommand{\supnorm}[2]{\ensuremath{\Vert{#1}\Vert}_{C(#2)}}

\begin{document}

\title[Density theorems with applications in QSP]{Density theorems with applications in\\ quantum signal processing}


\author{Rahul Sarkar}
\address{Institute for Computational and Mathematical Engineering, Stanford University, Stanford, CA 94305}
\curraddr{}
\email{rsarkar@stanford.edu}
\thanks{Rahul Sarkar was funded by the Stanford Exploration Project for the duration of this study.}

\author{Theodore J.~Yoder}
\address{IBM T.J. Watson Research Center, Yorktown Heights, NY}
\curraddr{}
\email{ted.yoder@ibm.com}
\thanks{}


\date{\today}

\dedicatory{}

\commby{}

\begin{abstract}
We study the approximation capabilities of two families of univariate polynomials that arise in applications of quantum signal processing. Although approximation only in the domain $[0,1]$ is physically desired, these polynomial families are defined by bound constraints not just in $[0,1]$, but also with additional bound constraints outside $[0,1]$. One might wonder then if these additional constraints inhibit their approximation properties within $[0,1]$. The main result of this paper is that this is not the case --- the additional constraints do not hinder the ability of these polynomial families to approximate arbitrarily well any continuous function $f:[0,1] \rightarrow [0,1]$ in the supremum norm, provided $f$ also matches any polynomial in the family at $0$ and $1$. We additionally study the specific problem of approximating the step function on $[0,1]$ (with the step from $0$ to $1$ occurring at $x=\frac{1}{2}$) using one of these families, and propose two subfamilies of monotone and non-monotone approximations. For the non-monotone case, under some additional assumptions, we provide an iterative heuristic algorithm that finds the optimal polynomial approximation.
\end{abstract}

\maketitle

\section{Introduction}
\label{sec:qsp-motivation}

In polynomial approximation theory, the famous Weierstrass approximation theorem \cite{weierstrass1885analytische} states that continuous functions over $G:=[0,1]$ can be uniformly approximated by polynomials in the supremum norm. As well as being a significant tool in functional analysis, the Weierstrass approximation theorem and its generalizations to cover the cases of generalized mappings due to Stone \cite{stone1948generalized}, called the Stone-Weierstrass theorem, and alternate topologies due to Krein \cite{krein1945problem}, are also practically useful, providing constructive approximations in signal processing \cite{carini2016study} and neural networks \cite{cotter1990stone}.

However, in some situations, one may be interested in approximations using a set of constrained polynomials in $G$, in which case the Weierstrass and Stone-Weierstrass theorems are inapplicable, and need to be generalized further. Some examples arise in stability analysis in control theory. For example, in the Lyapunov analysis of linear time-delay systems, one seeks polynomial solutions to optimization problems with affine constraints, and a generalization of the Weierstrass theorem in this setting to linear varieties can be found in \cite{peet2007extension}. Another recent work \cite{peet2009exponentially} used a generalization of the Weierstass theorem to Sobolev spaces to establish conditions under which an exponentially stable system has a polynomial Lyapunov function.

Instances of constrained polynomial approximation also arise in numerical analysis and in many scientific computing applications. Some examples of these are one-sided and comonotone polynomial approximations of a function \cite{trigub2009approximation,shvedov1980comonotone}, approximation by polynomials with positive coefficients \cite{trigub1998approximation,toland1996self}, and approximation by polynomials with integer coefficients \cite{le1980approximation,trigub1962approximation} (also see the references within \cite{trigub2009approximation} for an exhaustive list of works on these topics). However, the most well-known application relates to the situation when the polynomials are constrained to be positive in $G$, i.e. approximation by positive polynomials \cite{campos2019algorithms,allen2021bounds}. Such examples include positive cubic polynomial and spline interpolation \cite{butt1993preserving,schmidt1988positivity}, characterization of sums of squares \cite{lasserre2009moments,powers2011positive}, computer aided design with Bernstein and B\'{e}zier curves \cite{risler1992mathematical}, and non-negative approximation of hyperbolic equations \cite{toro2013riemann,zhang2013optimized}, to name a few. 

In all of the aforementioned cases, except in the cases of approximation by polynomials with positive or integer coefficients and comonotone approximations (which do not fall under the category of pointwise constraints), the approximating polynomials have pointwise constraints only on the set over which the approximation is desired, such as $G$. But what happens if one additionally has pointwise constraints outside of $G$? We study one such situation here, arising from quantum signal processing. Specifically, we consider the two classes of polynomials $\mathcal{P}$ and $\mathcal{Q}$, given in the following definition:

\begin{definition}
\label{def:polynomial-subset}
Let $\mathcal{P}$ be the set of real-valued, univariate polynomials $p$ satisfying the three additional properties: (i) $p(x)\in[0,1]$ for all $x\in G$, (ii) $p(x)\le0$ for $x\le0$, and (iii) $p(x)\ge1$ for $x\ge1$. Let $\mathcal{Q}$ be the set of real-valued, univariate polynomials $p$ satisfying properties (i) and (ii) of $\mathcal{P}$, and in addition $(\text{iii}^{\prime})$ $p(x)\le 0$ for $x\ge1$.
\end{definition}

What functions $f:[0,1]\rightarrow[0,1]$ can be approximated to arbitrarily small error in the supremum norm by polynomials in $\mathcal{P}$ or $\mathcal{Q}$? Just based on the behavior of the polynomials in sets $\mathcal{P}$ and $\mathcal{Q}$ over the domain $[0,1]$, it is at the very least necessary that $f(0)=0$, $f(x)\in[0,1]$ for $x\in[0,1]$, and $f(1)=1$ (or $f(1)=0$) for arbitrarily close approximation by $\mathcal{P}$ (or resp.~$\mathcal{Q}$). Our main result is that these conditions on $f$ are also sufficient. 

Some related work worth mentioning are \cite{hintermuller2015density,hintermuller2017density}, where density of convex subsets of some Banach spaces given by pointwise constraints have been studied. But the results derived in \cite{hintermuller2015density} are not directly applicable to our setting as the sets $\mathcal{P}$ and $\mathcal{Q}$ are not closed (they are only convex), while the pointwise constraints considered in \cite{hintermuller2017density} only apply to the set over which approximation is desired (which would be $G$ in our setting), and do not involve any constraints outside of $G$, such as those enforced by properties (ii), (iii), and $(\text{iii}^{\prime})$ in Definition~\ref{def:polynomial-subset}. Also of note is an algorithm \cite{campos2020projection} to project a function in $G$ to the set of polynomials satisfying only property (i) of Definition~\ref{def:polynomial-subset}, based on representation formulas in \cite{despres2017polynomials}.

To provide some context to our main result, let us briefly examine a set of polynomials $\tilde{\mathcal{P}}$ that are constrained exactly as $\mathcal{P}$ within the domain $[0,1]$ but differently outside, which results in $\tilde{\mathcal{P}}$ being much less capable of arbitrarily close approximation.
\begin{definition}
Let $a_0, a_1, b_0$, and $b_1$ be univariate polynomials satisfying $a_0(0)=a_1(0)=0$, $b_0(1)=b_1(1)=1$, $a_0(x)\le a_1(x)$ for all $x\le 0$, and $b_0(x)\le b_1(x)$ for all $x\ge1$. Let $\tilde{\mathcal{P}}$ be the set of univariate polynomials $p$ satisfying $p(x)\in[0,1]$ for all $x\in G$, $a_0(x)\le p(x)\le a_1(x)$ for all $x \le 0$, and $b_0(x)\le p(x)\le b_1(x)$ for all $x \ge 1$.
\end{definition}
Outside of $G$, polynomials in $\tilde{\mathcal{P}}$ are required to lie between polynomials $a_0$, $a_1$ for $x\le 0$, and between $b_0$, $b_1$ for $x\ge 1$. It is clear then that for any $p\in \tilde{\mathcal{P}}$, $\deg(p)\le d:= \max(\deg(a_0),\deg(a_1),\deg(b_0),\deg(b_1))$. However, it is known that for a non-polynomial function $f \in C(G)$, there exists $\epsilon_d>0$ such that no polynomial with degree at most $d$ can approximate $f$ over the interval $G$ to an error less than $\epsilon_d$ in the supremum norm (apply, for example, Theorem 3 in \cite{mayans2006chebyshev}). Thus no polynomial in $\tilde{\mathcal{P}}$ can approximate $f$ over $G$ to an error less than $\epsilon_d$ also (as $\tilde{\mathcal{P}}$ is a subset of all degree $d$ polynomials), and thus arbitrarily close approximation in the supremum norm to non-polynomial continuous functions is impossible using $\tilde{\mathcal{P}}$.

A common application of standard Chebyshev approximation theory is in filter design in digital signal processing, wherein some physically realizable class of functions, usually taken to be at least continuous, must be used to approximate ideal filter shapes (rectangular, triangular, etc). A similar problem in \emph{quantum signal processing} (QSP) motivates the study of $\mathcal{P}$ and $\mathcal{Q}$. In brief, in QSP one wishes to approximate a desired unitary function $U:\theta\in[0,\pi]\rightarrow\mathrm{SU}(2)$ by alternately applying a simple unitary function of $\theta$ and $\theta$-independent unitaries \cite{haah2019product,low2016methodology}. Next, we explain in more detail how the sets of polynomials $\mathcal{P}$ and $\mathcal{Q}$ arise from QSP.

\subsection{Constructibility questions in QSP}
\label{subsec:constructibility-qsp}

A quantum state of $n$ qubits is a unit-magnitude element from the Hilbert space $\mathbb{C}^{2^n}$ equipped with the standard Euclidean inner product, and operations on such a state are represented by unitary matrices (or operators) from $\mathrm{SU}(2^n)$. Some example unitaries from $\mathrm{SU}(2)$ are the Pauli matrices
\begin{equation}
I:=\left(\begin{array}{cc}1&0\\0&1\end{array}\right),\quad X:=\left(\begin{array}{cc}0&1\\1&0\end{array}\right),\quad Y:=\left(\begin{array}{cc}0&-i\\i&0\end{array}\right),\quad Z:=\left(\begin{array}{cc}1&0\\0&-1\end{array}\right),
\end{equation}
and these matrices are orthogonal with respect to the trace inner product, so they form a basis for $\mathbb{C}^{2 \times 2}$. 

In quantum signal processing, we deal with unitary functions $U:[0,\pi]\rightarrow \mathrm{SU}(2)$. Specifically, one is given access to a $1$-qubit \textit{signal} unitary function $S : \theta \mapsto \cos(\theta/2)I-i\sin(\theta/2)Y=e^{-i\theta Y/2} \in \mathrm{SU}(2)$ that encodes the signal $\theta\in[0,\pi]$, and asked to create a different 1-qubit unitary function $U$ such that for all values of $\theta \in [0,\pi]$, $U(\theta)$ can be formed using a finite product of $S(\theta)$ and some other prescribed set of $\theta$-independent $1$-qubit unitary operators. In our case, following \cite{low2016methodology}, we take the $\theta$-independent set to be a commutative multiplicative subgroup of $\mathrm{SU}(2)$
\begin{equation}
\mathcal{Z}:=\left\{Z(\phi):=\cos(\phi/2)I-i\sin(\phi/2)Z=e^{-i\phi Z/2}:\phi\in[0,2\pi)\right\}.
\end{equation}

Without loss of generality, because $\phi, \phi' \in [0, 2\pi)$ and $Z(0)=I$, a product of $L$ uses of $S(\theta)$ and elements of $\mathcal{Z}$ can be written as
\begin{equation}\label{eq:sequence}
V(\theta;\{\phi_0,\dots,\phi_L\})=Z(\phi_L)S(\theta)Z(\phi_{L-1})S(\theta)\dots Z(\phi_1)S(\theta)Z(\phi_0),
\end{equation}
where this product is parameterized by the angles $\phi_i\in[0,2\pi)$, for $i=0,1,2,\dots,L$. Any unitary function in a set
\begin{equation}
\mathcal{C}_L := \{U: [0,\pi] \ni \theta \mapsto V(\theta;\{\phi_0,\dots,\phi_L\}):\phi_0,\dots,\phi_L\in[0,2\pi)\}
\end{equation}
for some positive integer $L$ is called \emph{constructible}. We emphasize that for $U$ to be constructible, one must be able to use $S(\theta)$ and unitaries from $\mathcal{Z}$ to create $U(\theta)$ \emph{for all} values of $\theta\in[0,\pi]$ simultaneously --- it is a question of unitary functional equivalence.


If $U$ is a unitary function, it can be written as $U(\theta)=A(\theta)I+iB(\theta)X+iC(\theta)Y+iD(\theta)Z$, where $A$, $B$, $C$, and $D$ are real valued functions with domain $[0,\pi]$, uniquely determined by $U$. If $U$ is constructible, then $A$, $B$, $C$, and $D$ must have certain properties. For example, because $S(\theta)=e^{-i\theta/2}\frac12(I+Y)+e^{i\theta/2}\frac12(I-Y)$, we can deduce from \eqref{eq:sequence} that $A$, $B$, $C$, and $D$ are functions of $t_\theta:=e^{i\theta/2}$, and even more specifically that they are Laurent polynomials in $t_\theta$. That is, $A(\theta)=\tilde{A}(t_\theta):=\sum_{k=-L}^L a_k t_\theta^k$ for complex numbers $a_k$ depending only on the angles $\phi_0,\dots,\phi_L$ (and independent of $\theta$), where $\tilde{A}$ denotes the corresponding Laurent polynomial, and likewise for $B$, $C$, and $D$. The Laurent polynomials corresponding to $B$, $C$, and $D$ are denoted $\tilde{B}$, $\tilde{C}$ and $\tilde{D}$ respectively. Now the following theorem is known:
\begin{theorem}
\label{thm:constructibility}
Let $t_\theta := e^{-i\theta/2}$ for $\theta \in [0,\pi]$, and suppose we are given Laurent polynomials $\tilde{A}$, $\tilde{B}$, $\tilde{C}$, and $\tilde{D}$. Then $U: \theta \mapsto \tilde{A}(t_\theta)I + i \tilde{B}(t_\theta)X + i \tilde{C}(t_\theta)Y + i \tilde{D}(t_\theta)Z$ belongs to $\mathcal{C}_L$ 
if and only if the following conditions hold:
\begin{enumerate}[(i)]
    \item For all $z\in\mathbb{C}$ with $|z|=1$, $\tilde{A}(z)^2+\tilde{B}(z)^2+\tilde{C}(z)^2+\tilde{D}(z)^2=1$.
    \item The largest degree of these four Laurent polynomials is $L$.
    \item Each $\tilde{A},\tilde{B},\tilde{C},\tilde{D}$ is an even function if $L$ is even and an odd function if $L$ is odd. Being even or odd means that $\tilde{A}(z)=\tilde{A}(-z)$ and $\tilde{A}(z)=-\tilde{A}(-z)$ for all $z \in \mathbb{C}\setminus \{0\}$ respectively, and similarly for $\tilde{B}$, $\tilde{C}$, and $\tilde{D}$.
    \item $\tilde{A}$ and $\tilde{D}$ are reciprocal, i.e. $\tilde{A}(z)=\tilde{A}(1/z)$ and $\tilde{D}(z)=\tilde{D}(1/z)$, while $\tilde{B}$ and $\tilde{C}$ are anti-reciprocal, i.e. $\tilde{B}(z)=-\tilde{B}(1/z)$ and $\tilde{C}(z)=-\tilde{C}(1/z)$ for all $z \in \mathbb{C} \setminus \{0\}$.
\end{enumerate}
\end{theorem}
One can find the proof of this theorem in \cite{haah2019product}, and algorithms for finding the necessary parameters $\phi_0,\dots,\phi_L$ given a constructible $U(\theta)$ exist \cite{low2016methodology} and are efficient \cite{haah2019product,chao2020finding,dong2021efficient}.

Refinements of the constructibility question are important as well, because one might be interested only in some property of a unitary function $U$ and not its entirety. For example, starting with $U$ and the corresponding functions $A$ and $B$, consider the function $F: [0,\pi] \in \theta \mapsto B(\theta)^2+C(\theta)^2 \in \mathbb{R}$. This is physically relevant when $U(\theta)$ is applied to the quantum state $|0\rangle :=\left(\begin{smallmatrix}1\\0\end{smallmatrix}\right)$, which gives the state
\begin{equation}
U(\theta)|0\rangle = \left(\begin{array}{c}A(\theta)+iD(\theta)\\-C(\theta)+iB(\theta)\end{array}\right),
\end{equation}
and $F(\theta)=B(\theta)^2+C(\theta)^2$ is the absolute value squared of the Euclidean inner product of $U(\theta)|0\rangle$ with $|1\rangle :=\left(\begin{smallmatrix}0\\1\end{smallmatrix}\right)$, 
i.e. $F(\theta) = \left|\langle1|U(\theta)|0\rangle\right|^2$. Physically, $F(\theta)$ is the probability that a binary-valued measurement of the state $U(\theta)|0\rangle$ gives the final state $|1\rangle$. We say that $F$ is constructible with $L$ uses of $S$, if $F$ is obtained this way starting from some $U \in \mathcal{C}_L$.

Let $\lambda_\theta := \sin^2(\theta/2)$ and note $\theta\in[0,\pi]$ implies $\lambda_\theta \in [0,1]$, with the map $\theta \mapsto \lambda_\theta$ being bijective. It turns out that if $F$ is constructible with $L$ uses of $S$, then one can write $F(\theta) = p(\lambda_\theta)$, where $p$ is a polynomial of degree $L$ determined uniquely by $F$, with polynomial coefficients independent of $\theta$. Moreover $F$ is constructible with an odd number $L$ uses of $S$ if and only if $p \in \mathcal{P}$ and has degree $L$, while $F$ is constructible with an even number $L$ uses of $S$ if and only if $p \in \mathcal{Q}$ and has degree $L$. In context, properties (ii), (iii), and $(\text{iii}^\prime)$ of $\mathcal{P}$ and $\mathcal{Q}$ in Definition~\ref{def:polynomial-subset} may now appear strange as they involve values of $\lambda_\theta \not \in [0,1]$. However, they are natural consequences of the fact that $F(\theta)=B(\theta)^2+C(\theta)^2$, and can be derived using Theorem~\ref{thm:constructibility}. We refer to \cite{low2016methodology} for the proofs of these characterizations of $F$. 
Finally, an equivalent characterization of the sets $\mathcal{P}$ and $\mathcal{Q}$ in terms of polynomial sum of squares also exists, which we prove below.

\begin{lemma}
$p\in\mathcal{P}$ if and only if there exist odd polynomials $a,b,c,d$ such that $p(x)=c(\sqrt{x})^2+d(\sqrt{x})^2$ and $p(x)=1-a(\sqrt{1-x})^2-b(\sqrt{1-x})^2$.
\end{lemma}

\begin{proof}
The reverse implication is obvious, so we just need to show that $p\in\mathcal{P}$ implies these sum of squares decompositions. We factor
\begin{equation}
\label{eq:pfactorized}
p(x)=\alpha\prod_{i=0}^{n-1}(x-\beta_i)^{k_i}\prod_{j=0}^{m-1}\left((x-\gamma_j)^2+\delta_j^2\right)^{l_j},
\end{equation}
where $\alpha, \beta_i, \gamma_j, \delta_j$ are real for all $i,j$. Note that because $p(x)\in\mathcal{P}$, we can let $\beta_0=0$, which implies $k_0$ is odd and every other $k_i$ is even. Also, $\alpha>0$.

Now, repeated use of the identity
\begin{equation}
(r^2+s^2)(t^2+u^2)=(rt\pm su)^2+(ru\mp st)^2,
\end{equation}
which holds for any complex numbers $r,s,t,u$ and either choice of signs of the right hand side, reduces the product over $j$ in \eqref{eq:pfactorized} to a single sum of squares $\tilde c(x)^2+\tilde d(x)^2$, where $\tilde c$ and $\tilde d$ are polynomials in $x$. Thus, define
\begin{equation}
\begin{split}
c(\sqrt{x}) &= \sqrt{\alpha x^{k_0}}\left(\prod_{i=1}^{n-1}(x-\beta_i)^{k_i/2}\right)\tilde c(x), \\
\quad d(\sqrt{x}) &= \sqrt{\alpha x^{k_0}}\left(\prod_{i=1}^{n-1}(x-\beta_i)^{k_i/2}\right)\tilde d(x).
\end{split}
\end{equation}
These are odd polynomials and, by the prior discussion, $p(x)=c(\sqrt{x})^2+d(\sqrt{x})^2$.

Similarly, to find odd polynomials $a$ and $b$, let $\bar p(x)=1-p(1-x)\in\mathcal{P}$. Decompose $\bar p(x)=a(\sqrt{x})^2+b(\sqrt{x})^2$ as we did for $p$, and notice that it implies $p(x)=1-a(\sqrt{1-x})^2-b(\sqrt{1-x})^2$.
\end{proof}

\begin{lemma}
$q\in\mathcal{Q}$ if and only if there exist even polynomials $a,b,c,d$ such that $q(x)=x(1-x)\left(c(\sqrt{x})^2+d(\sqrt{x})^2\right)$ and $q(x)=1-a(\sqrt{x})^2-b(\sqrt{x})^2$.
\end{lemma}
\begin{proof}
Again, because the reverse implication is obvious, we just show that $q\in\mathcal{Q}$ implies the existence of $a,b,c,d$. We use the well-known polynomial sum-of-squares (SoS) theorem (see for instance \cite{marshall2008positive}): supposing $t:\mathbb{R}\rightarrow\mathbb{R}$ is a polynomial, $t(x)\ge0$ for all $x\in\mathbb{R}$ if and only if there are polynomials $u,v$ such that $t(x)=u(x)^2+v(x)^2$.

Note, because $q$ has roots at $x=0,1$, that $\tilde q(x)=q(x)/(x(1-x))$ is a polynomial and $\tilde q(x)\ge0$ for all $x\in\mathbb{R}$. The SoS theorem implies that $\tilde c$ and $\tilde d$ exist so that $q(x)=x(1-x)\left(\tilde c(x)^2+\tilde d(x)^2\right)$. Defining even polynomials $c$ and $d$ so that $c(\sqrt{x})=\tilde c(x)$ and $d(\sqrt{x})=\tilde d(x)$ is straightforward.

Likewise, $\bar q(x)=1-q(x)\ge 0$ for all $x\in\mathbb{R}$. The SoS theorem implies the existence of polynomials $\tilde a$ and $\tilde b$ such that $q(x)=1-\tilde a(x)^2-\tilde b(x)^2$, and even polynomials $a$ and $b$ exist such that $a(\sqrt{x})=\tilde a(x)$ and $b(\sqrt{x})=\tilde b(x)$.
\end{proof}


We note here that properties (ii) and (iii) in Definition~\ref{def:polynomial-subset} imply that all polynomials in $\mathcal{P}$ are of odd degree, while properties (ii) and $(\text{iii}^\prime)$ imply that all polynomials in $\mathcal{Q}$ are of even degree.

If one wants some $f: [0,1] \rightarrow \mathbb{R}$ that is \emph{not} in $\mathcal{P}$ or $\mathcal{Q}$ (so it cannot be constructed exactly), then one must resort to approximating $f$ (for instance, in the supremum norm) using the functions that are constructible, the polynomials in $\mathcal{P}$ and $\mathcal{Q}$. This motivates the question of which functions $f$ can be arbitrarily well approximated by $\mathcal{P}$ and $\mathcal{Q}$.

\subsection{Contributions}
\label{subsec:contributions}

Here we summarize the main contributions of this paper:
\begin{enumerate}[(i)]
    \item In Section~\ref{sec:density-result}, we prove the main result of our paper (Theorems~\ref{thm:densityP} and \ref{thm:densityQ}), which quantifies the closures of the polynomial classes $\mathcal{P}$ and $\mathcal{Q}$ in the supremum norm. These closures turn out to be the family of continuous functions $\mathcal{F}_1$ and $\mathcal{F}_2$ respectively (see Definition~\ref{def:cont-func-subset}), which are closed convex sets of the Banach space of continuous functions on $G$. The results extend easily to compact subsets of $G$, as in Corollaries~\ref{cor:density-closed-subsets-P} and \ref{cor:density-closed-subsets-Q}.
    \item In Section~\ref{sec:stepfunc-approx} we study the approximation of the step function $\Theta$ defined in \eqref{eq:step-func} over $A_\epsilon =[0,1/2-\epsilon]\cup[1/2+\epsilon,1]$, using elements of $\mathcal{P}$. It is shown in Section~\ref{ssec:monotone-approx} that Bernstein polynomial approximations are in $\mathcal{P}$ if and only if the polynomial degree is $1 \pmod{4}$, in which case it forms a monotonically increasing approximation to $\Theta$.
    \item In Section~\ref{ssec:nonmonotone-approx}, we propose a family $\mathcal{R}_\ell$ of non-monotone polynomial approximations to $\Theta$ over $A_\epsilon$ (Lemma~\ref{lem:type2-poly}), and then show that there exist elements of this family that are also elements of $\mathcal{P}$ (Lemma~\ref{lem:inP}). We further introduce a subfamily of $\mathcal{R}_\ell$ called equi-ripple polynomials (Definition~\ref{def:equiripple}), and  conjecture when they achieve the best approximation error to $\Theta$, among all the elements of $\mathcal{R}_\ell$. This is supported by numerical evidence, and we propose a heuristic algorithm (Algorithm~\ref{alg:iterative}) that we believe finds this best polynomial.
\end{enumerate}
\section{The Density Results}
\label{sec:density-result}

We first introduce some notation that will be useful. Recall $G:=[0,1]$. Let $H$ be a compact subset of $G$ and denote by $C(H)$ the Banach space of continuous real valued functions on $H$ equipped with the supremum norm, i.e. if $f \in C(H)$, then $\supnorm{f}{H} := \sup_{x \in H} |f(x)|$. 
Moreover, for ease of notation, if $G \subseteq H \subseteq \mathbb{R}$, $f \in C(G)$, and $g: H \rightarrow \mathbb{R}$, we define $f \pm g := f \pm (g|_G)$. 

We define $\mathcal{P}_H := \{p|_{H} : p \in \mathcal{P}\}$, i.e. the set of functions obtained by restricting each element of $\mathcal{P}$ to $H$. $\mathcal{Q}_H$ is similarly defined to be the restriction of elements of $\mathcal{Q}$ to $H$. With a slight abuse of notation, if $p \in \mathcal{P}$ or $\mathcal{Q}$, we will also denote its restriction to $H$ by $p$. Clearly $\mathcal{P}_H$, $\mathcal{Q}_H \subseteq C(H)$, and they are non-empty as the polynomial $p(x) = x$ belongs to $\mathcal{P}$, while $p(x)=x(1-x)$ belongs to $\mathcal{Q}$. We denote the closure of $\mathcal{P}_H$ and $\mathcal{Q}_H$ in $C(H)$ by $\overline{\mathcal{P}_H}^{C(H)}$ and $\overline{\mathcal{Q}_H}^{C(H)}$ respectively.

\begin{definition}
\label{def:cont-func-subset}
Let $\mathcal{F}_1$ be the set of all functions $f \in C(G)$ such that $f(0)=0$, $f(1)=1$, and $0 \leq f(x) \leq 1$ for all $x \in G$. Let $\mathcal{F}_2$ be the set of all functions $f \in C(G)$ such that $f(0)=f(1)=0$, and $0 \leq f(x) \leq 1$ for all $x \in G$.
\end{definition}

It is immediate from continuity that if $p \in \mathcal{P}$, then $p(0)=0$ and $p(1)=1$, while if if $p \in \mathcal{Q}$, then $p(0)=p(1)=0$. Next suppose $f,g \in \mathcal{F}_1$ (resp. $\mathcal{F}_2$) and let $\lambda \in [0,1]$. Then $\lambda f + (1-\lambda) g \in \mathcal{F}_1$ (resp. $\mathcal{F}_2$). Thus, the sets $\mathcal{F}_1$ and $\mathcal{F}_2$ are convex subsets of $C(G)$. Similarly $\mathcal{P}_G$ and $\mathcal{Q}_G$ are also convex subsets of $C(G)$. Moreover, $\mathcal{F}_1$ and $\mathcal{F}_2$ are closed in the topology\footnote{Here the topology is the metric space topology induced by the supremum norm.} of $C(G)$, because if we take a convergent sequence $f_i \rightarrow f \in C(G)$, where each $f_i \in \mathcal{F}_1$ (resp. $\mathcal{F}_2$), then we must have $0 \leq f(x) \leq 1$ for all $x \in G$, $f(0)=0$, and $f(1)=1$ (resp. $f(1)=0$), so in fact $f \in \mathcal{F}_1$ (resp. $\mathcal{F}_2$).

The main goal of this section is to establish the following theorems:

\begin{theorem}
\label{thm:densityP}
$\overline{\mathcal{P}_G}^{C(G)} = \mathcal{F}_1$.
\end{theorem}

\begin{theorem}
\label{thm:densityQ}
$\overline{\mathcal{Q}_G}^{C(G)} = \mathcal{F}_2$.
\end{theorem}

The proofs of these theorems are entirely constructive. Before we present the proofs, we collect some well-known results below that we need. The first of these is a constructive proof of the Weierstrass approximation theorem \cite{weierstrass1885analytische,perez2006survey} due to S. Bernstein \cite{bernstein1912}. The second result concerns a fact about piecewise monotone polynomial interpolation, which was proved independently by Young and Wolibner \cite{young1967piecewise,wolibner1951polynome}.

\begin{lemma}[Weierstrass approximation theorem]
\label{lem:bern-approx}
Let $f \in C(G)$, and for $n \geq 1$, define the degree-$n$ Bernstein polynomial approximation to $f$ by 
\begin{equation}
\label{eq:Bernstein-approx}
    (B_n f)(x) := \sum_{k=0}^{n} f\left( \frac{k}{n} \right) {n \choose k} x^k (1-x)^{n-k}, \;\; \text{ for all } x \in \mathbb{R}.
\end{equation}
Then the sequence $B_n f$ converges uniformly to $f$ on $G$, i.e. $\supnorm{B_n f - f}{G} \rightarrow 0$ as $n \rightarrow \infty$.
\end{lemma}

\begin{lemma}[Piecewise motonone polynomial interpolation]
\label{lem:piecewise-monotone-interp}
Suppose that $n$ is a positive integer, and we are given a set of points $\{(x_i,y_i) \in \mathbb{R}^2 : i=0,\dots,n\}$ satisfying $x_{i-1} < x_i$ and $y_{i-1} \neq y_i$ for $i=1,\dots,n$. Then there exists a polynomial $p$ such that $p(x_i)=y_i$, for $i=0,\dots,n$, and $p$ is monotone in each interval $[x_{i-1},x_i]$, for $i=1,\dots,n$.
\end{lemma}

We will also need the following lemma:
\begin{lemma}
\label{lem:correction-poly-step3}
Define the polynomial $r_{\alpha,\beta}(x) := x^{2\alpha + 1} (1-x)^{\beta}$, for all $\; x \in \mathbb{R}$, and integers $\alpha, \beta > 0$. Then $r_{\alpha,\beta}$ satisfies (a) $r_{\alpha,\beta}(0)=r_{\alpha,\beta}(1)=0$, $r_{\alpha,\beta}(x) > 0$ for $x \in (0,1)$, $r_{\alpha,\beta}(x) < 0$ for $x < 0$, and $r_{\alpha,\beta}(x) (-1)^\beta > 0$ for $x > 1$, (b) $r_{\alpha,\beta}$ is monotonically increasing for $x \leq 0$ and $r_{\alpha,\beta}(-1)^\beta$ is monotonically increasing for $x \geq 1$, and (c) $\sup_{x \in G} r_{\alpha,\beta}(x) = (1+\gamma)^{-2\alpha-1} (1 + \gamma^{-1})^{-\beta} < 2^{-\min \{2 \alpha+1, \beta\}}$, where $\gamma = \beta / (2\alpha + 1)$. Moreover,
\begin{enumerate}[(i)]
    \item if $\alpha_1, \alpha_2, \beta_1, \beta_2 > 0$ are positive integers with $\alpha_2 > \alpha_1$ and $\beta_2 > \beta_1$, then (a) $r_{\alpha_1,\beta_1}(x) \geq r_{\alpha_2,\beta_1}(x)$ and $r_{\alpha_1,\beta_1}(x) \geq r_{\alpha_1,\beta_2}(x)$ for all $x \in (0,1) \cup \{x \leq -1\}$, (b) $(r_{\alpha_1,\beta_1}(x) - r_{\alpha_2,\beta_1}(x)) (-1)^{\beta_1} \leq 0$ for all $x \geq 2$, and (c) $(-1)^{\beta_1} r_{\alpha_1,\beta_1}(x) \leq  (-1)^{\beta_2} r_{\alpha_1,\beta_2}(x)$ for all $x \geq 2$.
    \item if $h$ is a polynomial, there exist positive integers $\alpha,\beta$ (with $\beta$ even) such that $r_{\alpha,\beta}(x) + h(x) > 0$ for all $x \geq 2$.
    \item if $h$ is a polynomial, there exist positive integers $\alpha,\beta$ (with $\beta$ odd) such that $r_{\alpha,\beta}(x) + h(x) < 0$ for all $x \geq 2$.
    \item if $h$ is a polynomial, there exist positive integers $\alpha,\beta$ such that $r_{\alpha,\beta}(x) + h(x) < 0$ for all $x \leq -1$.
\end{enumerate}
\end{lemma}

\begin{proof}
We first prove the first part, so let $\alpha,\beta > 0$ be integers. From the definition of $r_{\alpha,\beta}$ it is clear that (a) holds. Differentiating $r_{\alpha,\beta}$ gives 
\begin{equation}
\label{eq:r-alpha-beta-prime}
\begin{split}
    r'_{\alpha,\beta}(x) &= (2\alpha + 1) x^{2 \alpha} (1-x)^{\beta} - \beta x^{2 \alpha + 1} (1-x)^{\beta - 1} \\
    &= x^{2\alpha} (1-x)^{\beta-1}((2\alpha + 1) - (\beta + 2 \alpha + 1) x),    
\end{split}
\end{equation}
and so the only critical points of $r_{\alpha,\beta}$ are $x=0,1,(1+\gamma)^{-1}$. Noticing that $(1+\gamma)^{-1} \in (0,1)$, we get from \eqref{eq:r-alpha-beta-prime} that $r'_{\alpha,\beta}(x) > 0$ for all $x < 0$, and $r'_{\alpha,\beta}(x) (-1)^\beta > 0$ for all $x > 1$. This proves (b). To prove (c), note that using (a) one concludes $r_{\alpha,\beta}$ attains its maximum on $G$ at an interior point, which must also be a critical point as $r_{\alpha,\beta}$ is smooth. But the only such point is $(1+\gamma)^{-1}$, and thus $\sup_{x \in G} r_{\alpha,\beta}(x) = r_{\alpha,\beta} ((1+\gamma)^{-1})=(1+\gamma)^{-2\alpha-1} (1 + \gamma^{-1})^{-\beta}$. Now there are two cases: (i) $\gamma > 1$ which gives $(1+\gamma)^{-1} < 1/2$ and $(1+\gamma^{-1})^{-1} < 1$, and (ii) $0 < \gamma \le 1$ in which case $(1+\gamma)^{-1} < 1$ and $0 < (1+\gamma^{-1})^{-1} \leq 1/2$, and in both these cases the bound $\sup_{x \in G} r_{\alpha,\beta}(x) < 2^{-\min \{2 \alpha+1, \beta\}}$ follows.

We now prove the remaining parts of the lemma.

(i) For $x \in (0,1)$, we have $0 < x^{2\alpha_2+1} < x^{2\alpha_1+1}$, and $0 < (1-x)^{\beta_2} < (1-x)^{\beta_1}$. For $x \leq -1$, we have $(1-x)^{\beta_2} > (1-x)^{\beta_1} \geq 2$, and $x^{2\alpha_2+1} \leq x^{2\alpha_1+1} \leq -1$. For $x \geq 2$, we have $(1-x)^{\beta_2} (-1)^{\beta_2} \geq (1-x)^{\beta_1} (-1)^{\beta_1}$, and $x^{2\alpha_2+1} > x^{2\alpha_1+1} \geq 8$. Combining these facts implies (a), (b), and (c).

(ii), (iii) Let $h$ be a polynomial of degree $n$, and $u \geq 1$. Then $h(1+u)= \sum_{i=0}^{n} c_i u^i$, for some constants $c_0,\dots,c_n$ independent of $u$, and this implies $|h(1+u)| \leq \sum_{i=0}^{n} |c_i| u^i \leq \left(\sum_{i=0}^{n} |c_i| \right) u^n = K u^n$, where $K = \sum_{i=0}^{n} |c_i|$. Now choose $\beta > n$ such that $\beta$ is even (resp. odd) for part (ii) (resp. (iii)), and $\alpha$ such that $2^{2\alpha+1} > K$. This gives $r_{\alpha,\beta}(1+u) (-1)^\beta = (1+u)^{2\alpha+1} u^{ \beta} \geq 2^{2\alpha+1} u^{\beta} > K u^n$. Thus $(r_{\alpha,\beta}(1+u) + h(1+u)) (-1)^\beta \geq (-1)^\beta r_{\alpha,\beta}(1+u) - |h(1+u)| \geq (-1)^\beta r_{\alpha,\beta}(1+u) - K u^n > 0$.

(iv) Let $h(x) = \sum_{i=0}^{n} c_i x^i$ be the polynomial, and fix $x \leq -1$. Then $|h(x)| \leq \sum_{i=0}^{n} |c_i x^i| \leq K |x|^n$, where where $K = \sum_{i=0}^{n} |c_i|$. Now choose $\alpha > n/2$, and $\beta$ such that $2^{\beta} > K$. This implies $r_{\alpha,\beta}(x) = x^{2\alpha+1} (1-x)^{\beta} \leq -|x|^n (1-x)^{\beta} \leq -|x|^n 2^{\beta} < -K|x|^n < 0$. Using these we deduce $r_{\alpha,\beta}(x) + h(x) \leq r_{\alpha,\beta}(x) + |h(x)| \leq r_{\alpha,\beta}(x) + K|x|^n < 0$.
\end{proof}

We now return to the proofs of Theorem~\ref{thm:densityP} and Theorem~\ref{thm:densityQ}. The main work involved in proving Theorem~\ref{thm:densityP} is in showing that $\mathcal{P}_G$ is dense in $\mathcal{F}_1$, in the topology of $C(G)$ (Theorem~\ref{thm:density1-P}), i.e. given any $f \in \mathcal{F}_1$ there exists a sequence $\{p_i \in \mathcal{P}_G : i=1,2,\dots\}$ such that $\supnorm{p_i - f}{G} \rightarrow 0$ as $i \rightarrow \infty$. Similarly, we also need to show that $\mathcal{Q}_G$ is dense in $\mathcal{F}_2$ (Theorem~\ref{thm:density1-Q}) for proving Theorem~\ref{thm:densityQ}.

\subsection{Proof of Theorem~\ref{thm:densityP}}
\label{ssec:densityP}

We will break up the proof of Theorem~\ref{thm:densityP} into several steps. The first step is given any $f \in \mathcal{F}_1$, we want to replace it by an arbitrarily close polynomial approximation $p$ satisfying certain properties, as stated in the next lemma.

\begin{lemma}
\label{lem:bern-approx-step1-P}
Let $f \in \mathcal{F}_1$, and $\epsilon > 0$. Then there exists a polynomial $p$ satisfying (i) $\supnorm{p - f}{G} < \epsilon$, (ii) $0 < p(x) < 1$ for all $x \in (0,1)$, (iii) $p(0)=0$, and (iv) $p(1)=1$.
\end{lemma}

\begin{proof}
First let $n \geq 1$ be arbitrary, and consider the Bernstein polynomial approximation $B_n f$ to $f$ as defined in Lemma~\ref{lem:bern-approx}. Then clearly $(B_n f)(0)=f(0)$, and $(B_n f)(1)=f(1)$. We also have $\sum_{k=0}^{n} {n \choose k} x^k (1-x)^{n-k} = 1$ for any $x \in G$, which implies $0=\inf_{x \in G} f(x) \leq (B_n f)(x) \leq \sup_{x \in G} f(x)=1$. Moreover by Lemma~\ref{lem:bern-approx}, there exists $N \geq 1$, such that $\supnorm{B_n f - f}{G} < \epsilon$ for every $n > N$. Now consider the set
\begin{equation}
    S := \bigcup_{n > N} \{k/n : k = 0,\dots,n\},
\end{equation}
and we claim that there exists $y \in S$ such that $0 < f(y) < 1$. Let us first show how the claim implies the lemma. Choose $n > N$ and  $k \in \{0,\dots,n\}$ such that $k/n \in S$ with $0 < f(k/n) < 1$. Then for all $x \in (0,1)$, we have $(B_n f)(x) = \sum_{k=0}^{n} f\left( \frac{k}{n} \right) {n \choose k} x^k (1-x)^{n-k} < \sum_{k=0}^{n} {n \choose k} x^k (1-x)^{n-k}=1$, and similarly we also have $(B_n f)(x) > 0$. Thus $B_n f$ satisfies all properties (i)-(iv) and is the required polynomial $p$.

To prove the claim, for contradiction suppose it is false. Then one can write $S = S_0 \sqcup S_1$, with the property that $f(y)=0$ for all $y \in S_0$, and $f(y)=1$ for all $y \in S_1$. Since $S$ is dense in $G$, at least one of the two sets $S_0$ or $S_1$ must also be dense in $G$. If $S_0$ (resp. $S_1$) is dense in $G$, the continuity of $f$ implies $f=0$ (resp. $f=1$) identically on $G$, which is a contradiction. 
\end{proof}

The polynomial $p$ furnished by Lemma~\ref{lem:bern-approx-step1-P} does not necessarily belong to $\mathcal{P}$. To correct for this, our next goal is to replace $p$ by another arbitrarily close approximating polynomial $q$, with correct properties over the interval $[-1,2]$. The precise result is stated in the next lemma, which uses monotone polynomial interpolation.

\begin{lemma}
\label{lem:monotone-approx-step2-P}
Let $\epsilon > 0$, and $p$ be a polynomial satisfying $p(0)=0$, $p(1)=1$, and $0 < p(x) < 1$ for all $x \in (0,1)$. Then there exists a polynomial $q$ satisfying the properties: (i) $\supnorm{q-p}{G} < \epsilon$, (ii) $0 < q(x) < 1$ for all $x \in (0,1)$, (iii) $q(0)=0$, (iv) $q(1)=1$, (v) $q(x) < 0$ for all $x \in [-1,0)$, and (vi) $q(x) > 1$ for all $x \in (1,2]$.
\end{lemma}

\begin{proof}
Let $p_G$ be the restriction of $p$ to $G$. Since $p_G \in C(G)$ and $G$ is compact, $p_G$ is uniformly continuous on $G$. Thus there exists $\delta > 0$ such that for all $x,y \in G$ with $|x-y| < \delta$, we have $|p_G(x) - p_G(y)| < \epsilon / 2$. We first claim that there exists a finite set of points $S := \{x_i : i=1,\dots,N\}$ satisfying (i) $0=x_1 < x_2 < \dots < x_{N-1} < x_N=1$, (ii) $x_{i} - x_{i-1} < \delta$ for all $i=2,\dots,N$, and (iii) $p_G(x_{i-1}) \neq p_G(x_{i})$ for all $i=2,\dots,N$. To prove the claim, we prescribe an iterative process. We start with $S = \{x_1: x_1=0\}$ and iteratively add a new element $x_{|S|+1}$ to $S$ while $1 - x_{|S|} \geq \delta$. The new element is chosen so that $\delta/2 \leq x_{|S|+1} - x_{|S|} < \delta$, and $p_G(x_{|S|}) \neq p_G(x_{|S|+1})$. Such a choice is possible because $p_G$ cannot be equal to $p_G(x_{|S|})$ everywhere on $[x_{|S|} + \delta/2, x_{|S|} + \delta)$, because otherwise $p_G$ would be a constant polynomial which cannot be in $\mathcal{F}_1$. The iterative process must terminate, as at any stage of this process we have $x_{|S|} \geq (|S|-1) \delta/2$. Once the process terminates, we simply add the point $x_N=1$ to $S$. Note that since $x_{N-1} \in (0,1)$, we have $p_G(x_{N-1}) \neq p_G(x_N)$. The set $S$ thus satisfies properties (i)-(iii) of the claim, and the claim is proved.

Now take any set $S = \{x_1,\dots,x_N\}$ as in the claim above, and define $x_0 = -1$ and $x_{N+1}=2$. Define $y_0=-1$, $y_{N+1}=2$, and $y_i = p_G(x_i)$ for all $i=1,\dots,N$. By Lemma~\ref{lem:piecewise-monotone-interp}, there exists a polynomial $q$ such that $q(x_i)=y_i$ for all $i=0,\dots,N+1$, and $q$ is monotone in each interval $[x_{i-1},x_i]$ for all $i=1,\dots,N+1$. We want to show that $q$ satisfies properties (i)-(vi) of the lemma. Clearly properties (iii) and (iv) are true by construction. Properties (v) and (vi) are consequences of monotonicity of $q$ on the intervals $[x_0,x_1]$ and $[x_N,x_{N+1}]$ respectively. Property (ii) again follows from monotonicity of $q$ in the remaining intervals, since $0 < y_i < 1$ for all $i=2,\dots,N-1$. To see that $q$ also satisfies property (i), we first fix $x \in G$ and suppose $x \in [x_{i-1},x_i]$ for some $i \in \{2,N\}$. Then we have 
\begin{equation}
\begin{split}
    |q(x) - p_G(x)| &= |q(x) - p_G(x_{i-1}) + p_G(x_{i-1}) - p_G(x)| \\
    & \leq |q(x) - p_G(x_{i-1}) | + |p_G(x_{i-1}) - p_G(x)| \\
    &= |q(x) - q(x_{i-1}) | + |p_G(x_{i-1}) - p_G(x)| \\
    & \leq |q(x_i) - q(x_{i-1}) | + |p_G(x_{i-1}) - p_G(x)| \\
    &= |p_G(x_i) - p_G(x_{i-1}) | + |p_G(x_{i-1}) - p_G(x)| \\
    &< \epsilon/2 + \epsilon /2 = \epsilon,
\end{split}
\end{equation}
where the first inequality is the triangle inequality, the second inequality is true by monotonicity of $q$, and the third inequality is due to uniform continuity of $p_G$ on $G$. Since $x \in G$ is arbitrary, the lemma is proved.
\end{proof}

Lemma~\ref{lem:bern-approx-step1-P} and Lemma ~\ref{lem:monotone-approx-step2-P} now allow us to prove the following theorem.

\begin{theorem}
\label{thm:density1-P}
$\mathcal{P}_G$ is dense in $\mathcal{F}_1$ in the topology of $C(G)$.
\end{theorem}

\begin{proof}
Let $f \in \mathcal{F}_1$, $\delta > 0$, and  $\epsilon := \delta / 3$. We want to show that there exists an element $u \in \mathcal{P}_G$ such that $\supnorm{u-f}{G} < \delta$. First we find a polynomial $p$ satisfying the properties in Lemma~\ref{lem:bern-approx-step1-P}, and then using this $p$, we find another polynomial $q$ satisfying the properties in Lemma~\ref{lem:monotone-approx-step2-P}. Thus at this stage we have $\supnorm{q-f}{G} \leq \supnorm{q-p}{G} + \supnorm{p-f}{G} < 2\epsilon$. If $q \in \mathcal{P}$ at this stage, then we are done. If not, we claim that there is a polynomial $r_{\alpha,\beta}$ as defined in Lemma~\ref{lem:correction-poly-step3}, such that $q + r_{\alpha,\beta} \in \mathcal{P}$, and $\sup_{x \in G} |r_{\alpha,\beta}(x)| < \epsilon$. If we define $u := q + r_{\alpha,\beta}$, then $\supnorm{u-f}{G} \leq \supnorm{q-f}{G} + \sup_{x \in G} |r_{\alpha,\beta}(x)| < 3\epsilon = \delta$, and the theorem is proved. It remains to prove the claim, which is done next.

Since $q$ is piecewise monotonic (with finitely many pieces in $[-1,2]$), $q(1)=1$, and $q((0,1)) \subseteq (0,1)$, there exists $\mu \in (0,1)$ such that $q$ is monotonically increasing in $[1-\mu,1]$. 
Let $H:=[0,1-\mu]$, and define $\eta := \sup_{x \in H} q(x)$. 
Since $q(H) \subseteq [0,1)$ and $H$ is compact, we have $q(1-\mu) \leq \eta < 1$.
Let $\eta' := \max (\eta, 1-\epsilon)$. 
Then by monotonicity of $q$ in $[1-\mu,1]$, there exists $\mu' \in [1-\mu,1)$ such that $q(\mu')=(1+\eta')/2 < 1$. 
Thus we have found $\mu' \in (0,1)$ with the properties: (i) $q(x) < q(\mu')$ for all $x \in [0,\mu')$, (ii) $q$ is monotonically increasing in $[\mu',1]$, and (iii) $1-q(\mu') < \epsilon$. 
Consider the polynomial $s := 1-q$ on the set $[\mu',1]$. We know $s(1)=0$ and $s(x) > 0$ for all $x \in [\mu',1)$. Thus one can factor $s$ as $s(x)=(1-x)^m \tilde{s}(x)$, where $\tilde{s}(1) \neq 0$, and $\tilde{s}(x) > 0$ for all $x \in [\mu',1)$. By continuity of $\tilde{s}$ and since $[\mu',1]$ is compact, there exists a lower bound $0 < K < \tilde{s}(x)$ for all $x \in [\mu',1]$, which implies $s(x) \geq K(1-x)^m$. Now choose $\beta_0$ such that $(1-\mu')^{\beta_0} < K$. Then for all $x \in [\mu',1)$, we have $(1-x)^{m+\beta_0} < K(1-x)^m$, and so $q(x) + (1-x)^{m+\beta_0} < q(x) + s(x) = 1$. 
Next, by property (c) of the first part of Lemma~\ref{lem:correction-poly-step3}, there exists a polynomial $r_{\alpha_1,\beta_1}$ such that $\sup_{x \in G} r_{\alpha_1,\beta_1}(x) < 1 - q(\mu')$.
Finally, if we set $h = q - q(2)$ in Lemma~\ref{lem:correction-poly-step3}(ii), we conclude that there exists a polynomial $r_{\alpha_2,\beta_2}$ such that $r_{\alpha_2,\beta_2}(x) + q(x) > q(2)=2$ for all $x \geq 2$, and similarly if we set $h = q - q(-1)$ in Lemma~\ref{lem:correction-poly-step3}(iv), we conclude that there exists a polynomial $r_{\alpha_3,\beta_3}$ satisfying $r_{\alpha_3,\beta_3}(x) + q(x) < q(-1)=-1$ for all $x \leq -1$. Moreover we can choose $\beta_0$, $\beta_1$, $\beta_2$ and $\beta_3$ to be even.

Define $\alpha := \sum_{i=1}^3 \alpha_i$, and $\beta := 2m + \sum_{i=0}^3 \beta_i$. By Lemma~\ref{lem:correction-poly-step3}(i), $\sup_{x \in G} r_{\alpha,\beta}(x) \leq \sup_{x \in G} r_{\alpha_1,\beta_1}(x)$, so $\sup_{x \in G} |r_{\alpha,\beta}(x)| = \sup_{x \in G} r_{\alpha,\beta}(x) < (1-q(\mu')) < \epsilon$ by property (iii) of $\mu'$ above, so we only need to show $q + r_{\alpha,\beta} \in \mathcal{P}$. Clearly by properties of $r_{\alpha,\beta}$ in the first part of Lemma~\ref{lem:correction-poly-step3}, we have $q(0)+r_{\alpha,\beta}(0)=0$, $q(1)+r_{\alpha,\beta}(1)=1$, and $q(x) + r_{\alpha,\beta}(x) > 0$ for all $x \in (0,1)$. Next for all $x \in [0,\mu']$, we have $q(x)+r_{\alpha,\beta}(x) \leq q(\mu') + r_{\alpha_1,\beta_1}(x) \leq q(\mu') + \sup_{x \in G} r_{\alpha_1,\beta_1}(x) < 1$, where the first inequality is by property (i) of $\mu'$ and Lemma~\ref{lem:correction-poly-step3}(i), while for all $x \in [\mu',1)$, we get $q(x)+r_{\alpha,\beta}(x) < q(x) + (1-x)^{m+\beta_0} \leq 1$ (by the last paragraph). Thus we have proved $(q + r_{\alpha,\beta})((0,1)) \subseteq (0,1)$. From monotonicity of $q$, and again using the properties of $r_{\alpha,\beta}$, we next obtain $q(x) + r_{\alpha,\beta}(x) > q(1) = 1$ for all $x \in (1,2]$, and $q(x) + r_{\alpha,\beta}(x) < q(0) = 0$ for all $x \in [-1,0)$. Finally, another application of Lemma~\ref{lem:correction-poly-step3}(i) gives $q(x) + r_{\alpha,\beta}(x) \geq q(x) + r_{\alpha_2,\beta_2}(x) > 2$ for all $x \geq 2$, and $q(x) + r_{\alpha,\beta}(x) \leq q(x) + r_{\alpha_3,\beta_3}(x) < -1$ for all $x \leq -1$. This proves the claim.
\end{proof}

We now finish the proof of Theorem~\ref{thm:densityP}.
\begin{proof}[Proof of Theorem~\ref{thm:densityP}]
We first show that $\overline{\mathcal{P}_G}^{C(G)} \subseteq \mathcal{F}_1$. Take a sequence $p_i \rightarrow f \in C(G)$, where each $p_i \in \mathcal{P}_G$. Since for all $i$, we have $p_i(0)=0$, $p_i(1)=1$, and $0 \leq p_i(x) \leq 1$ for all $x \in G$, this implies $f(0)=0$, $f(1)=1$, and $0 \leq f(x) \leq 1$ for all $x \in G$, because otherwise it contradicts $\supnorm{p_i - f}{G} \xrightarrow{i \uparrow \infty} 0$. This proves $f \in \mathcal{F}_1$, and so $\overline{\mathcal{P}_G}^{C(G)} \subseteq \mathcal{F}_1$. Furthermore, by Theorem~\ref{thm:density1-P} we also have that $\mathcal{P}_G$ is dense in $\mathcal{F}_1$, so $\mathcal{F}_1 \subseteq \overline{\mathcal{P}_G}^{C(G)}$. The theorem is proved.
\end{proof}

If we study the proof of Theorem~\ref{thm:density1-P} carefully, the properties of $q + r_{\alpha,\beta}$ actually implies that we have also proved the following theorem:

\begin{theorem}
\label{thm:density2-P}
Let $\mathcal{P}'$ be the set of polynomials $p$ satisfying the properties: (i) $p(x)\in (0,1)$ for all $x\in (0,1)$, (ii) $p(x) < 0$ for all $x < 0$, and (iii) $p(x) > 1$ for all $x > 1$. Define $\mathcal{P}'_G := \{p|_{G} : p \in \mathcal{P}'\}$. Then $\mathcal{P}'_G$ is dense in $\mathcal{F}_1$ in the topology of $C(G)$, and $\overline{\mathcal{P}'_G}^{C(G)} = \mathcal{F}_1$.
\end{theorem}

We state an easy corollary of Theorem \ref{thm:density2-P} below.

\begin{corollary}
\label{cor:density-closed-subsets-P}
Let $A$ be a closed subset of $G$, $\epsilon > 0$, and $\mathcal{P}'$ be defined as in Theorem~\ref{thm:density2-P}. Let $f \in C(A)$ such that $0 \leq f(x) \leq 1$ for all $x \in A$, $f(0)=0$ if $0 \in A$, and $f(1)=1$ if $1 \in A$. Then there exists $p \in \mathcal{P}'$ with $\sup_{x \in A} |p(x) - f(x)| < \epsilon$.
\end{corollary}

\begin{proof}
As $G$ is a normal topological space, by the Tietze extension theorem (Theorem 35.1 in \cite{munkrestopology}) there exists $g \in \mathcal{F}_1$ that extends $f$. By Theorem \ref{thm:density2-P} there exists $p \in \mathcal{P}'$ with $\sup_{x \in A} |p(x) - g(x)| < \epsilon$.
\end{proof}

\subsection{Proof of Theorem~\ref{thm:densityQ}}
\label{ssec:densityQ}

The proof of Theorem~\ref{thm:densityQ} follows along similar lines. The first two steps are analogs of Lemma~\ref{lem:bern-approx-step1-P} and Lemma~\ref{lem:monotone-approx-step2-P}.

\begin{lemma}
\label{lem:bern-approx-step1-Q}
Let $f \in \mathcal{F}_2$, $f(x) \neq 0$ for some $x \in G$, and $\epsilon > 0$. Then there exists a polynomial $p$ satisfying (i) $\supnorm{p - f}{G} < \epsilon$, (ii) $0 < p(x) < 1$ for all $x \in (0,1)$, (iii) $p(0)=0$, and (iv) $p(1)=0$.
\end{lemma}

\begin{proof}
The exact same proof of Lemma~\ref{lem:bern-approx-step1-P} works as we have assumed $f$ is not identically zero on $G$.
\end{proof}

\begin{lemma}
\label{lem:monotone-approx-step2-Q}
Let $\epsilon > 0$, and $p$ be a non-constant polynomial satisfying $p(0)=p(1)=0$, and $0 < p(x) < 1$ for all $x \in (0,1)$. Then there exists a polynomial $q$ satisfying the properties: (i) $\supnorm{q-p}{G} < \epsilon$, (ii) $0 < q(x) < 1$ for all $x \in (0,1)$, (iii) $q(0)=0$, (iv) $q(1)=0$, (v) $q(x) < 0$ for all $x \in [-1,0)$, and (vi) $q(x) < 0$ for all $x \in (1,2]$.
\end{lemma}

\begin{proof}
The proof is the same as Lemma~\ref{lem:monotone-approx-step2-P} with the only changes that we use $\mathcal{F}_2$ in place of $\mathcal{F}_1$ and set $y_{N+1} = - 1$. The existence of the set $S$ also uses the assumption that $p$ is not a constant polynomial.
\end{proof}

We can now prove:
\begin{theorem}
\label{thm:density1-Q}
$\mathcal{Q}_G$ is dense in $\mathcal{F}_2$ in the topology of $C(G)$.
\end{theorem}

\begin{proof}
This proof is very similar to the proof of Theorem~\ref{thm:density1-P}, but there are some differences. Let $f \in \mathcal{F}_2$, $\delta > 0$, and  $\epsilon := \delta / 3$. If $f = 0$ identically on $G$, take $u := \delta x(1-x)$. Then we have $u \in \mathcal{Q}$ and $\supnorm{u-f}{G} < \delta$. So now assume  in the rest of the proof that $f$ is not identically zero on $G$. First we find polynomials $p$ and $q$ using Lemma~\ref{lem:bern-approx-step1-Q} and Lemma~\ref{lem:monotone-approx-step2-Q} so that $\supnorm{q-f}{G} < 2 \epsilon$. Then if $q \not \in \mathcal{Q}$, we claim that there exists $r_{\alpha,\beta}$ such that $u := q + r_{\alpha,\beta} \in \mathcal{Q}$ and $\sup_{x \in G} |r_{\alpha,\beta}(x)| < \epsilon$. Then $\supnorm{u-f}{G} < 3 \epsilon$, proving the theorem.

For the claim define $\eta := \sup_{x \in G} q(x)$, and let $\eta' := \max (\eta, 1- \epsilon)$. Since $G$ is compact, $q(0)=q(1)=0$, and $q((0,1)) \subseteq (0,1)$, we have $0 < \eta, \eta' < 1$. By property (c) of the first part of Lemma~\ref{lem:correction-poly-step3}, there exists $r_{\alpha_1,\beta_1}$ such that $\sup_{x \in G} r_{\alpha_1,\beta_1}(x) < 1- \eta'$. Setting $h = q - q(2)$ in Lemma~\ref{lem:correction-poly-step3}(iii), we obtain a polynomial $r_{\alpha_2,\beta_2}$ such that $r_{\alpha_2,\beta_2}(x) + q(x) < q(2)=-1$ for all $x \geq 2$, and setting $h = q - q(-1)$ in Lemma~\ref{lem:correction-poly-step3}(iv), we obtain a polynomial $r_{\alpha_3,\beta_3}$ satisfying $r_{\alpha_3,\beta_3}(x) + q(x) < q(-1)=-1$ for all $x \leq -1$. Moreover we can take $\beta_1$, $\beta_2$, and $\beta_3$ to be odd. Now define $\alpha := \sum_{i=1}^3 \alpha_i$, and $\beta := \sum_{i=1}^3 \beta_i$. Then $\sup_{x \in G} |r_{\alpha,\beta}(x)| \leq \sup_{x \in G} r_{\alpha_1,\beta_1}(x) < 1- \eta' \leq \epsilon$, where the first inequality is by Lemma~\ref{lem:correction-poly-step3}(i). To show that $q + r_{\alpha,\beta} \in \mathcal{Q}$, we note that the first part of Lemma~\ref{lem:correction-poly-step3} gives $q(0)+r_{\alpha,\beta}(0)=q(1)+r_{\alpha,\beta}(1)=0$, and $q(x) + r_{\alpha,\beta}(x) > 0$ for all $x \in (0,1)$. Next from the definitions of $\eta, \eta'$, for $x \in G$ we have $q(x) + r_{\alpha,\beta}(x) \leq \sup_{x \in G} q(x) + \sup_{x \in G} r_{\alpha_1,\beta_1}(x) < \eta + 1- \eta' < 1$. Finally, by properties (v)-(vi) of $q$ in Lemma~\ref{lem:monotone-approx-step2-Q}, property (a) of $r_{\alpha,\beta}$ in the first part of Lemma~\ref{lem:correction-poly-step3} (noting that $\beta$ is odd), and Lemma~\ref{lem:correction-poly-step3}(i), we have  that $q(x) + r_{\alpha,\beta}(x) < 0$ for all $x \in [-1,0) \cup (1,2]$, and $q(x) + r_{\alpha,\beta}(x) < -1$ for all $x \leq -1$ and $x \geq 2$.
\end{proof}

The proof of Theorem~\ref{thm:densityQ} now follows exactly similarly as that of Theorem~\ref{thm:densityP}:
\begin{proof}[Proof of Theorem~\ref{thm:densityQ}]
Theorem~\ref{thm:density1-Q} proves the inclusion $\mathcal{F}_2 \subseteq \overline{\mathcal{Q}_G}^{C(G)}$. For the opposite inclusion $\overline{\mathcal{Q}_G}^{C(G)} \subseteq \mathcal{F}_2$, we take a convergent sequence $p_i \rightarrow f \in C(G)$ where each $p_i \in \mathcal{Q}_G$. Then it must be true that $f \in \mathcal{F}_2$ due to the properties of each $p_i$, since convergence is in the supremum norm.
\end{proof}

We also immediately deduce the following theorem and corollary from the proof of Theorem~\ref{thm:densityQ}.

\begin{theorem}
\label{thm:density2-Q}
Let $\mathcal{Q}'$ be the set of polynomials $p$ satisfying the properties: (i) $p(x)\in (0,1)$ for all $x\in (0,1)$, (ii) $p(x) < 0$ for all $x < 0$, and (iii) $p(x) < 0$ for all $x > 1$. Define $\mathcal{Q}'_G := \{p|_{G} : p \in \mathcal{Q}'\}$. Then $\mathcal{Q}'_G$ is dense in $\mathcal{F}_2$ in the topology of $C(G)$, and $\overline{\mathcal{Q}'_G}^{C(G)} = \mathcal{F}_2$.
\end{theorem}

\begin{corollary}
\label{cor:density-closed-subsets-Q}
Let $A$ be a closed subset of $G$, $\epsilon > 0$, and $\mathcal{Q}'$ be defined as in Theorem~\ref{thm:density2-Q}. Let $f \in C(A)$ such that $0 \leq f(x) \leq 1$ for all $x \in A$, $f(0)=0$ if $0 \in A$, and $f(1)=0$ if $1 \in A$. Then there exists $p \in \mathcal{Q}'$ with $\sup_{x \in A} |p(x) - f(x)| < \epsilon$.
\end{corollary}

\begin{proof}
Replace $\mathcal{F}_1$ with $\mathcal{F}_2$, $\mathcal{P}'$ with $\mathcal{Q}'$, and use Theorem~\ref{thm:density2-Q} instead of Theorem~\ref{thm:density2-P} in the proof of Corollary~\ref{cor:density-closed-subsets-P}.
\end{proof}

\section{Approximating the step function}
\label{sec:stepfunc-approx}

For some $0 < \epsilon < 1/2$, consider the closed subset $A_\epsilon := [0,1/2-\epsilon]\cup[1/2+\epsilon,1]$, and define the \textit{step-function}
\begin{equation}
\label{eq:step-func}
\Theta(x):=
\begin{cases}
    0,&x<1/2 \\
    \frac{1}{2},&x=1/2 \\
    1,&x>1/2
\end{cases}
\end{equation}
on $G$. The step function is important in quantum signal processing as it answers the question: is the signal, here $x$, less than 1/2 or greater? In contexts where one wants to learn the signal by extracting one bit of information at a time, such as Hamiltonian simulation, and phase or amplitude estimation \cite{low2017hamiltonian,rall2021faster}, constructing the step function is a core subroutine.

We know already from Corollary~\ref{cor:density-closed-subsets-P} that for any given $\delta > 0$, there exists $p \in \mathcal{P}$ such that $\sup_{x \in A_\epsilon} |p(x) - \Theta (x)| < \delta$. In fact starting from the function $\Theta|_{A_\epsilon}$, the restriction of $\Theta$ to $A_\epsilon$, we can extend it to a continuous function $f \in \mathcal{F}_1$, and then the proof of Theorem~\ref{thm:density1-P} gives us a constructive way of obtaining such a $p$. In this section we furnish two other (and more direct) ways of finding polynomial approximations from $\mathcal{P}$ to the function $\Theta$ over the domain $A_\epsilon$. In Section~\ref{ssec:monotone-approx} we introduce a class of monotone polynomial approximations to this problem based on the Bernstein polynomial approximation \eqref{eq:Bernstein-approx}. In Section~\ref{ssec:nonmonotone-approx} we discuss a class of non-monotone polynomial approximations, and based on some unproven assumptions about this class, we also propose a heuristic algorithm for computing these polynomials.

\subsection{Monotone approximations}
\label{ssec:monotone-approx}

A simple way of constructing a monotone approximation to $\Theta$ by polynomials in $\mathcal{P}$ is to consider the degree-$L$ Bernstein polynomial approximation

\begin{equation}
\label{eq:Bernstein-step-func}
    (B_L \Theta)(x) := \sum_{k=0}^{L} \Theta \left( \frac{k}{L} \right) {L \choose k} x^k (1-x)^{L-k}, \;\; \text{ for all } x \in \mathbb{R}.
\end{equation}
For $B_L \Theta$ to be an element of $\mathcal{P}$, we already know that $L$ has to be odd. In fact one can say more, as we prove in the next lemma.

\begin{lemma}
\label{lem:bernstein-pmod4}
Let $B_L \Theta$ be defined as in \eqref{eq:Bernstein-step-func} and suppose $L$ is an odd positive integer. Then we have the following:
\begin{enumerate}[(i)]
    \item $B_L \Theta$ is monotonically increasing in $G$.
    \item $(B_L \Theta)(x) + (B_L \Theta) (1-x) = 1$, for all $x \in \mathbb{R}$.
    \item $B_L \Theta \in \mathcal{P}$ if and only if $L \equiv 1 \pmod{4}$.
\end{enumerate}
\end{lemma}

\begin{proof}
(i) The arguments in \cite[Section 2] {gzyl2003approximation} prove that $0 = (B_L \Theta) (0) \leq (B_L \Theta) (x) \leq (B_L \Theta) (1) = 1$ for $x \in G$, and that $B_L \Theta$ is non-decreasing in $G$. Since $B_L \Theta$ is a polynomial, we can also conclude that there does not exist distinct $x,y \in G$ such that $(B_L \Theta) (x) = (B_L \Theta) (y)$, as the non-decreasing property would otherwise imply that $B_L \Theta$ is a constant. Thus $B_L \Theta$ must be strictly increasing in $G$.

(ii) Let $L = 2\ell + 1$ for some integer $\ell$. Then by \eqref{eq:Bernstein-step-func} we get for $x \in \mathbb{R}$
\begin{equation*}
\begin{split}
    & (B_L \Theta)(x) + (B_L \Theta)(1-x) = \sum_{k=\ell + 1}^{2\ell + 1} {2\ell+1 \choose k} \left[ x^k (1-x)^{2\ell+1 - k} + (1-x)^k x^{2\ell+1 - k} \right] \\
    &= \sum_{k=\ell + 1}^{2\ell + 1} {2\ell+1 \choose k} x^k (1-x)^{2\ell+1 - k} + \sum_{k=0}^{\ell} {2\ell+1 \choose k} x^k (1-x)^{2\ell+1 - k} = 1,
\end{split}
\end{equation*}
using the binomial theorem. 

(iii) From the proof of part (i) we already have $0 = (B_L \Theta) (0) \leq (B_L \Theta) (x) \leq (B_L \Theta) (1) = 1$ for all $x \in G$. We will prove that $(B_L \Theta)(x) \leq 0$ for all $x < 0$, if and only if $L \equiv 1 \pmod{4}$. Part (ii) then implies $(B_L \Theta)(x) \geq 1$ for all $x > 1$, if and only if $L \equiv 1 \pmod{4}$, which would prove the result.
First assume $L = 4 \ell + 1$ for some integer $\ell$, and fix $x < 0$. If $\ell = 0$, then $(B_L \Theta)(x) = x < 0$, so we may assume $\ell \geq 1$. Then \eqref{eq:Bernstein-step-func} gives
\begin{equation}
\label{eq:bern-step-1}
\begin{split}
    (B_L \Theta) (x) &= x^{4 \ell + 1} + \sum_{k=0}^{\ell - 1} {4 \ell + 1 \choose 2\ell + 2k + 1} x^{2\ell + 2k + 1} (1-x)^{2\ell - 2k} \\
    &+ \sum_{k=0}^{\ell - 1} {4 \ell + 1 \choose 2\ell + 2k + 2} x^{2\ell + 2k + 2} (1-x)^{2 \ell - 2k - 1} \\
    &< \sum_{k=0}^{\ell - 1} \frac{(4 \ell + 1)! \; x^{2\ell + 2k + 1} (1-x)^{2 \ell - 2k - 1} (2\ell + 2k + 2 -x (4k+2)) }{(2\ell + 2k + 2)! \; (2\ell - 2k)!} \\
    &< 0.
\end{split}
\end{equation}
Next assume $L = 4 \ell + 3$ for some integer $\ell$, and fix $x < 0$. Then again using \eqref{eq:Bernstein-step-func} we get
\begin{equation}
\label{eq:bern-step-2}
\begin{split}
    (B_L \Theta) (x) &= \sum_{k=0}^{\ell} \frac{(4 \ell + 3)! \; x^{2\ell + 2k + 2} (1-x)^{2\ell - 2k} (2\ell + 2k + 3 -x (4k + 2))}{(2\ell + 2k + 3)! \; (2\ell - 2k + 1)!} \\
    & > 0,
\end{split}
\end{equation}
which finishes the proof.
\end{proof}

Using the probabilistic interpretation of the Bernstein polynomial approximation \cite[Section 2]{gzyl2003approximation}, and invoking results from large deviations theory, one can precisely state how fast $B_L \Theta$ converges to $\Theta$ on $A_\epsilon$:

\begin{lemma}
\label{lem:bernstein-step-rate}
Let $B_L \Theta$ be defined as in \eqref{eq:Bernstein-step-func} and suppose $L \equiv 1 \pmod{4}$. Then $\supnorm{B_L \Theta - \Theta}{A_\epsilon} \leq 2 e^{-2L \epsilon^2}$.
\end{lemma}

\begin{proof}
By Lemma~\ref{lem:bernstein-pmod4}(i),(ii) we have $\supnorm{B_L \Theta - \Theta}{A_\epsilon} = (B_L \Theta) (\frac{1}{2} - \epsilon)$. It then follows from the last paragraph of \cite[Section 2]{gzyl2003approximation} that $(B_L \Theta) (\frac{1}{2} - \epsilon) \leq 2 e^{-2L \epsilon^2}$.
\end{proof}

\subsection{Non-monotone approximations}
\label{ssec:nonmonotone-approx}

Notice that for $x \in A_\epsilon$, the function $\Theta$ exhibits the symmetry $\Theta(x)=1-\Theta(1-x)$. Moreover the monotone approximations $B_L \Theta$, for $L$ odd and $L \equiv 1 \pmod{4}$ that we introduced in Section~\ref{ssec:monotone-approx}, also exhibits the same symmetry (Lemma~\ref{lem:bernstein-pmod4}(ii)). Thus while seeking non-monotone polynomial approximations to $\Theta$ in $\mathcal{P}$, it is sensible to start from the class of polynomials with this symmetry. With this as the motivation, the following lemma provides the basic building block of our construction:

\begin{lemma}
\label{lem:type2-poly}
Suppose we have real numbers $0 < a_1 < a_2 < \dots < a_\ell < \frac{1}{2}$, for some integer $\ell \geq 1$. Then there exists a unique polynomial $p$ of degree $L := 4\ell + 1$, that satisfies
\begin{enumerate}[(i)]
    \item $p(x) + p(1-x) = 1$, for all $x \in \mathbb{R}$,
    \item $p(0) = p(a_i) = p'(a_i)=0$ for $i=1,2,\dots,\ell$, where $p'$ is the derivative of $p$.
\end{enumerate}

\vspace{0.2cm}
\noindent
Moreover, the polynomial $p$ also satisfies
\begin{enumerate}[(i)]
    \item[(iii)] $p'$ has exactly one zero in $(0,a_1)$, and in $(a_i,a_{i+1})$ for each $i=1,\dots,\ell-1$, and $p'(x) > 0$ for all $x \in (a_\ell, \frac{1}{2}]$ and $x \leq 0$, 
    \item[(iv)] $p(x) < 0$ for all $x < 0$, and $p(x) \geq 0$ for all $x \in [0,\frac{1}{2}]$ with equality if and only if $x=0$ or $x=a_i$ for some $i$.
\end{enumerate}
\end{lemma}

\begin{proof}
Define $\bar{a}_i := 1-a_i$ for $i=1,2,\dots,\ell$. Then we have $\frac{1}{2} < \bar{a}_\ell < \bar{a}_{\ell-1} < \dots < \bar{a}_2 < \bar{a}_1 < 1$. For uniqueness note that if a polynomial $p$ satisfies properties (i),(ii), then we deduce that $p$ must also satisfy the conditions: $\text{(ii)}^\prime$ $p(1) = p(\bar{a}_i) = 1$ for $i=1,\dots,\ell$, and $\text{(ii)}^{\prime \prime}$ $p'(\bar{a}_i)=0$ for $i=1,2,\dots,\ell$. Since $L = 4\ell + 1$, the discussion in \cite[Section 5]{spitzbart1960generalization} shows that $p$ is unique. For existence, we note that by \cite[Theorem 1 and Section 5]{spitzbart1960generalization}, there exists a unique polynomial $p$ of degree $4\ell + 1$ satisfying properties (ii), $\text{(ii)}^\prime$, and $\text{(ii)}^{\prime \prime}$. We will show that $p$ satisfies (i). For all $x \in \mathbb{R}$, let $p(x) = \sum_{k=0}^{4 \ell + 1} c_k (x - \frac{1}{2})^k$, for some real-valued coefficients $c_0,\dots,c_{4\ell+1}$. Define $a_0 := 0$ and $\bar{a}_0 := 1$. Then properties (ii), $\text{(ii)}^\prime$, and $\text{(ii)}^{\prime \prime}$ lead to the system of linear equations
\begin{equation}
\label{eq:nonmonotone-system1}
\begin{split}
    \sum_{k=0}^{4 \ell + 1} c_k \left( z - 1/2 \right)^k = 0, \; \sum_{k=0}^{4 \ell + 1} c_k \left( \bar{z} - 1/2 \right)^k = 1, \; &\text{ for } z=a_0,\dots,a_{\ell}, \\
    \sum_{k=1}^{4 \ell+1} k c_k \left( z - 1/2 \right)^{k-1} = \sum_{k=1}^{4 \ell+1} k c_k \left( \bar{z} - 1/2 \right)^{k-1} = 0, \; & \text{ for } z=a_1,\dots,a_{\ell}.
\end{split}
\end{equation}
Eliminating the coefficients $\{c_k : k \text{ odd}\}$ from \eqref{eq:nonmonotone-system1} leads to the linear system
\begin{equation}
\label{eq:nonmonotone-system2}
\begin{split}
    \sum_{\substack{k=0 \\ k \text{ even}}}^{4 \ell + 1} c_k \left( z - 1/2 \right)^k = 1/2, \; & \text{ for } z=a_0,\dots,a_{\ell},   \\
    \sum_{\substack{k=1 \\ k \text{ even}}}^{4 \ell+1} k c_k \left( z - 1/2 \right)^{k-1} = 0, \; & \text{ for } z=a_1,\dots,a_{\ell}.
\end{split}
\end{equation}
Now by uniqueness of $p$, \eqref{eq:nonmonotone-system1} has a unique solution for the coefficients $c_0,\dots,c_{4\ell+1}$, and hence the solution to \eqref{eq:nonmonotone-system2} is unique too. Noticing that $c_0 = \frac{1}{2}$, and $c_k = 0$ for $k = 2,4,\dots,4\ell$ satisfies \eqref{eq:nonmonotone-system2}, we conclude that $p$ must be of the form $p(x) = \frac{1}{2} + \sum_{k=0}^{2 \ell} c_{2k+1} (x - \frac{1}{2})^{2k+1}$, and it is then obvious that $p$ satisfies (i).

We now show that $p$ satisfies properties (iii), (iv). By Rolle's theorem, property (ii) implies that for each $i=1,\dots,\ell$, there exists $b_i \in (a_{i-1},a_i)$ such that $p'(b_i)=0$. Then property (i) gives $p'(\bar{b}_i)=0$ for $i=1,\dots,\ell$, where $\bar{b}_i = 1-b_i$, and $\bar{b}_i \in (\bar{a}_{i},\bar{a}_{i-1})$. Combining with $\text{(ii)}^{\prime \prime}$ we find that $p'$ has $4\ell$ distinct real zeros in $(0,1)$, and since $p'$ is a degree-$4\ell$ polynomial, these are all its zeros, each of multiplicity one. This implies $p(x) = 0$ for $x \leq \frac{1}{2}$ if and only if $x \in \{a_i: i=0,\dots,\ell\}$, because otherwise $p'$ will have at least $4\ell+1$ distinct real zeros, which is not possible. Now fix some $a_i$, for $i=1,\dots,\ell$. Then $a_i$ is a zero of $p$ of multiplicity at least two due to property (ii), but the multiplicity also cannot exceed two, as $a_i$ is a zero of $p'$ of multiplicity one. Thus either $p(x) \geq 0$ or $p(x) \leq 0$ for all $x \in [0,\frac{1}{2}]$, but since property (i) implies $p(\frac{1}{2})=\frac{1}{2}$, we conclude the former. Next, $p'$ must be strictly positive or negative for all $x < b_1$, as $b_1$ is the smallest zero of $p'$. Assuming it is strictly negative, it leads to a contradiction because by the previous paragraph $p(b_1) > 0$, and $p(0)=0$ by property (ii). Hence we now conclude $p(x) < 0$ for all $x < 0$. Finally note that $p'$ has no zero in $(a_\ell,\bar{a}_\ell)$, and hence must be strictly positive or strictly negative in this open interval. But if it is strictly negative, then it contradicts $p(\frac{1}{2})=\frac{1}{2}$.
\end{proof}

Any polynomial $p$ obtained from the above lemma satisfies properties (ii) and (iii) of the polynomial set $\mathcal{P}$ in Definition~\ref{def:polynomial-subset}, but it does not necessarily satisfy property (i) of the definition, because it may be that $p(x) > 1$ for some $x \in (0,\frac{1}{2})$. So the question is if one can choose the real numbers $0 < a_1 < a_2 < \dots < a_\ell < \frac{1}{2}$ in such a way, such that $p(x) \leq 1$, for all $x \in [0,\frac{1}{2}]$. If this can be done, then $p \in \mathcal{P}$ due to property (i) of Lemma~\ref{lem:type2-poly}. The next lemma shows that this is indeed possible:

\begin{lemma}
\label{lem:inP}
There exist real numbers $0 < a_1 < a_2 < \dots < a_\ell < \frac{1}{2}$ for which the polynomial $p$ of degree $4\ell + 1$ satisfying properties (i) and (ii) of Lemma~\ref{lem:type2-poly}, also satisfies $p(x) \leq 1$ for all $x \in [0,\frac{1}{2}]$.
\end{lemma}

\begin{proof}
Let $a_0 := 0$ and $\bar{a}_0 := 1$. As in the proof of Lemma~\ref{lem:type2-poly}, for each $i=1,\dots,\ell$, let $b_i \in (a_{i-1},a_{i})$ be such that $p'(b_i)=0$, where $p'$ denotes the derivative of $p$, and also define $\bar{b}_i := 1-b_i \in (\bar{a}_i, \bar{a}_{i-1})$ and $\bar{a}_i := 1-a_i$. Then we know from the proof that the set $\bigcup_{i=1}^{\ell} \{a_i, \bar{a}_i, b_{i}, \bar{b}_{i}\}$ are all the zeros of $p'$. So let $p'$ have the form
\begin{equation}\label{eq:pprime_factorization}
    p'(x) = A \prod_{i=1}^{\ell} (x- a_i) (x - \bar{a}_i) (x- b_i) (x - \bar{b}_i), \;\; x \in \mathbb{R},
\end{equation}
and $A > 0$, since $p'(x) > 0$ for all $x \leq 0$ by Lemma~\ref{lem:type2-poly}(iii). Notice that properties (i) and (ii) of $p$ in Lemma~\ref{lem:type2-poly} imply that for each $i=1,\dots,\ell$, we have $\int_{a_{i-1}}^{a_i} p'(x) \; dx = \int_{\bar{a}_{i}}^{\bar{a}_{i-1}} p'(x) \; dx = 0$, and since $p(0)=0$ and $p(1)=1$, we may further conclude that
\begin{equation}
\label{eq:pprime-expressions1}
    p(y) = \int_{0}^{y} p'(x) \; dx = \int_{a_i}^{y} p'(x) \; dx, \;\; \text{ if } y \in [a_i,a_{i+1}], \; i=0,\dots,\ell-1,
\end{equation}
and
\begin{equation}
\label{eq:pprime-expressions2}
    A \int_{a_{\ell}}^{\bar{a}_{\ell}} \prod_{i=1}^{\ell} (x- a_i) (x - \bar{a}_i) (x- b_i) (x - \bar{b}_i) \; dx = \int_{a_{\ell}}^{\bar{a}_{\ell}} p'(x) \; dx = p(1) = 1.
\end{equation}
Now choose $a_\ell < w < \frac{1}{2}$, and define $\bar{w}:= 1-w$. Then we have the following estimate:
\begin{equation}
\label{eq:estimate1-p}
\begin{split}
    & \int_{a_{\ell}}^{\bar{a}_{\ell}} \prod_{i=1}^{\ell} (x- a_i) (x - \bar{a}_i) (x- b_i) (x - \bar{b}_i) \; dx \\
    > & \int_{w}^{\bar{w}} \prod_{i=1}^{\ell} |x- a_i| \; |x - \bar{a}_i| \; |x- b_i| \; |x - \bar{b}_i| \; dx \\
    > & \int_{w}^{\bar{w}} \prod_{i=1}^{\ell} (w - a_\ell)^4 \; dx = (w - a_\ell)^{4 \ell} (1 - 2 w).
\end{split}
\end{equation}
Combining \eqref{eq:estimate1-p} and \eqref{eq:pprime-expressions2} gives $A < (w - a_\ell)^{-4 \ell} (1 - 2 w)^{-1}$. Next pick any $j \in \{1,\dots,\ell\}$. Then using \eqref{eq:pprime-expressions1}, this bound on $A$ gives the estimate
\begin{equation}
\label{eq:estimate2-p}
\begin{split}
    p(b_j) &= A \int_{a_{j-1}}^{b_j} \prod_{i=1}^{\ell} (x- a_i) (x - \bar{a}_i) (x- b_i) (x - \bar{b}_i) \; dx \\
    & \leq A \int_{a_{j-1}}^{b_j} \prod_{i=1}^{\ell} |x- a_i| \; |x - \bar{a}_i| \; |x- b_i| \; |x - \bar{b}_i| \; dx \\
    & \leq A \; (b_j - a_{j-1}) <  (w - a_\ell)^{-4 \ell} (1 - 2 w)^{-1} a_\ell,
\end{split}
\end{equation}
where the third inequality is true because each of the terms $|x- a_i|$, $|x - \bar{a}_i|$, $|x- b_i|$, and $|x - \bar{b}_i|$ are smaller than $1$. Now keeping $w$ fixed, if we let $a_\ell \rightarrow 0$, then by \eqref{eq:estimate2-p} we get that $p(b_j) \rightarrow 0$. Since $j$ is arbitrary, and because we proved in Lemma~\ref{lem:type2-poly} that $b_j$ is the unique maximum of $p$ in $[a_{j-1},a_j]$, we can conclude that there exists $a_1,\dots,a_\ell$ such that $p(x) \leq 1$ for all $x \in [0,a_\ell]$. Finally, for $x \in (a_\ell,\frac{1}{2}]$, we know that $p'(x) > 0$ by Lemma~\ref{lem:type2-poly}(iii), and so $p$ is monotonically increasing in $[a_\ell,\frac{1}{2}]$, which implies $p(x) \leq \frac{1}{2}$.
\end{proof}

From this point onward, we specialize further to the class of \emph{equi-ripple} polynomials. These are defined as follows. First, denote polynomials obtained using Lemma~\ref{lem:type2-poly} for fixed $\ell \geq 1$ by $\mathcal{R}_\ell$.

\begin{definition}\label{def:equiripple}
Suppose $a_1,\dots,a_\ell$ be as in Lemma~\ref{lem:type2-poly}, and let $a_0 := 0$. Let $p \in \mathcal{R}_\ell$, and for every $i = 1,\dots,\ell$, let $b_i$ be the unique point where $p$ achieves its maximum value in $[a_{i-1},a_{i}]$. We then say that $p$ has the \textit{equi-ripple property} if for some $\delta > 0$, we have $p(b_i)=p(b_j)=\delta$ for all $i,j \in \{1,\dots,\ell\}$. We often refer to $\delta$ as the \textit{ripple height}.
\end{definition}

Polynomials with the equi-ripple property are similar to best approximating polynomials in standard Chebyshev approximation theory. The difference here is that the error is one-sided — an equi-ripple polynomial is lower bounded by $\Theta$ in $[0,\frac{1}{2}]$. 

We claim that the polynomial $p$ in $\mathcal{R}_\ell$ best approximates $\Theta$ in the supremum norm over $A_\epsilon$, meaning that $\supnorm{p - \Theta}{A_\epsilon} \leq \supnorm{q - \Theta}{A_\epsilon}$ for every $q \in \mathcal{R}_\ell$, if (i) $p$ has the equi-ripple property, (ii) $a_\ell < \frac{1}{2} - \epsilon$, and (iii) $p(\frac{1}{2}-\epsilon)=\delta$. Note, this automatically implies $\delta < \frac{1}{2}$ and thus $p\in\mathcal{P}$, since the polynomial is monotonically increasing from $\frac{1}{2}-\epsilon$ to $\frac{1}{2}$, and $p(\frac{1}{2})=\frac{1}{2}$. A question now emerges regarding the existence of such a polynomial with properties (i), (ii), and (iii). We claim that if one fixes $\ell$ (resp. fixes $a_\ell$), there is a sufficiently small $a_\ell$ (resp. a sufficiently large $\ell$) such that all three properties hold. These claims are backed-up by numerical studies. We have implemented an iterative algorithm that we believe always terminates with an equi-ripple polynomial.

The idea of this iterative algorithm is as follows. Let $a_0:= 0$. Starting with fixed $\ell \geq 1$ and $a_\ell$, we adjust the locations of $a_i$ for $i=1,2,\dots,\ell-1$ until the maxima (located at $b_i$ in Definition~\ref{def:equiripple}) are all of the same height. More specifically, start with any initial values of $a_1,\dots,a_{\ell-1} \in (0,a_\ell)$. In each iteration, we find the current largest maximum (at $b_M$ say) and current smallest maximum (at $b_m$ say). If the largest and smallest maxima are of the same height, then we have achieved the equi-ripple property and the algorithm can terminate. Otherwise, widen the interval $I_m=[a_{m-1},a_m]$ at the expense of the interval $I_M=[a_{M-1},a_M]$, while keeping all other intervals the same width, until the maxima within these intervals $I_m$ and $I_M$ are of the same height. That this can always be done is guaranteed by the following lemma, which establishes upper and lower bounds on the height of a peak $p(b_i)$ in terms of the width of the interval $I_i$.

\begin{lemma}
If $b_i$ is the unique point in $[a_{i-1},a_i]$ where a polynomial $p\in\mathcal{R}_\ell$ achieves its maximum value, then
\begin{equation}\label{eq:squeeze_peak}
\left[4(4\ell+1)\binom{4\ell}{2i-1}\right]^{-1}\frac{(a_i-a_{i-1})^{4\ell-1}}{1/2-a_\ell}<p(b_i)<2^{4\ell}\frac{a_i-a_{i-1}}{(1/2-a_\ell)^{4\ell+1}}.
\end{equation}
\end{lemma}

\begin{proof}
We get by \eqref{eq:estimate2-p} from the proof of Lemma~\ref{lem:inP} that $p(b_i)\le A(b_i-a_{i-1})$, and so $p(b_i)<A(a_i-a_{i-1})$ also. Moreover, from that proof, we see $A$ is upper bounded by a constant --- specifically $A<(w-a_\ell)^{-4\ell}(1-2w)^{-1}$, where $w$ is any number satisfying $a_\ell<w<1/2$. Choosing $w=(1/2+a_\ell)/2$ gives the upper bound in \eqref{eq:squeeze_peak}.

To get the lower bound, we look at two cases, (i) $b_i\ge (a_i+a_{i-1})/2$ and (ii) $b_i<(a_i+a_{i-1})/2$. In case (i), we start with \eqref{eq:pprime_factorization} and calculate for $x \in [a_{i-1}, b_i]$,
\begin{equation}
\begin{split}
p'(x)&=A\prod_{j=1}^\ell(x-a_j)(x-\bar a_j)(x-b_j)(x-\bar b_j)\\
&=A\prod_{j=1}^{i-1}|x-a_j||x-b_j|\prod_{j=i}^\ell|x-a_j||x-b_j|\prod_{j=1}^\ell|x-\bar a_j||x-\bar b_j|\\
&\ge A|x-a_{i-1}|^{2(i-1)}|x-b_i|^{2(\ell-i+1)}|x-b_i|^{2\ell}.
\end{split}
\end{equation}
We therefore find
\begin{equation}
\label{eq:casei_LB}
\begin{split}
p(b_i)=\int_{a_{i-1}}^{b_i}p'(x)dx&\ge A\int_{a_{i-1}}^{b_i}(x-a_{i-1})^{2i-2}(x-b_i)^{4\ell-2i+2}dx\\
&=A\left[(4\ell+1)\binom{4\ell}{2i-2}\right]^{-1}(b_i-a_{i-1})^{4\ell+1}\\
&\ge A\left[(4\ell+1)\binom{4\ell}{2i-2}\right]^{-1}\left(\frac{a_i-a_{i-1}}{2}\right)^{4\ell+1}.
\end{split}
\end{equation}

A similar calculation handles case (ii). We find for $x \in [b_i,a_i]$
\begin{align}
p'(x)&\ge A|x-b_i|^{2i-1}|x-a_i|^{2\ell-2i+1}|x-a_i|^{2\ell},\\\label{eq:caseii_LB}
p(b_i)& \geq A\left[(4\ell+1)\binom{4\ell}{2i-1}\right]^{-1}\left(\frac{a_i-a_{i-1}}{2}\right)^{4\ell+1}.
\end{align}
Note, since $A>0$, the lower bound from case (ii) in \eqref{eq:caseii_LB}, is always smaller than that from case (i) in \eqref{eq:casei_LB}. The final piece is a lower bound on $A$. We use \eqref{eq:pprime-expressions2} and the bound
\begin{equation}
\prod_{i=1}^\ell(x-a_i)(x-\bar a_i)(x-b_i)(x-\bar b_i)<1/2^{4\ell}
\end{equation}
to get $A>2^{4\ell}/(1-2a_\ell)$, which completes the proof.
\end{proof}

After adjusting the intervals $I_m$ and $I_M$, the algorithm repeats, finding new intervals containing the largest and smallest peaks. See the pseudo-code in Algorithm~\ref{alg:iterative}. 
In practice, the algorithm has some smallest step size $\kappa$ by which it is willing to adjust the intervals. We set $\kappa=10^{-4}$ in our experiments.
We observe that Algorithm~\ref{alg:iterative} successfully returns an equi-ripple polynomial for some range of parameters $\ell$ and $a_{\ell}$. When $\ell$ is too large or $a_{\ell}$ is too close to $\frac{1}{2}$, the algorithm can struggle to find a polynomial with zeros at the specified $a_i$ (line 3 in Algorithm~\ref{alg:iterative}). In principle, this step involves solving a system of linear equations for the coefficients of the polynomial $p$. However in practice, the matrix inversion can be badly conditioned. Extending the practicality of Algorithm~\ref{alg:iterative} is an interesting open question.

\begin{algorithm}
\caption{Return a set of zeros (with additive error $\kappa>0$) so that ripples are roughly equal in height.}\label{alg:remez}
\begin{algorithmic}[1]
    \Procedure{Iterate}{$a=\{a_1,a_2,\dots,a_\ell\}$,$\kappa$}
        \State Let $a_0=0$
        \State Find polynomial $p\in\mathcal{R}_\ell$ with zeros $a_i$ (polynomial from Lemma~\ref{lem:type2-poly})
        \State Find the heights $h_1,h_2\dots,h_{\ell}$ of local maxima of $p$ in $(0,a_\ell)$
        \State If $h_M$ is the largest, let $b_M$ be the peak's location
        \State If $h_m$ is the smallest, let $b_m$ be the peak's location
        \State Set $s=\min((a_m-a_{m-1})/4,(a_M-a_{M-1})/4)$
        \State If $b_M<b_m$, set $s=-s$
        \State While $|s|\ge\kappa$
        \State \quad For $j=1,2,\dots,\ell$
        \State \quad \quad If $\min(m,M)\le j< \max(m,M)$, let $a'_j=a_j+s$
        \State \quad \quad Else $a'_j=a_j$
        \State \quad Find polynomial $p\in\mathcal{R}_\ell$ with zeros $a'_i$ (polynomial from Lemma~\ref{lem:type2-poly})
        \State \quad Find the heights $h_1,h_2\dots,h_{\ell}$ of local maxima of $p$ in $(0,a_\ell)$
        \State \quad If $h_m\le h_M$, set $a_j=a'_j$ for all $j$
        \State \quad If $h_m>h_M$
        \State \quad \quad If $|s|=\kappa$, set $s=0$
        \State \quad \quad Else, set $s=\text{sign}(s)\max(|s/2|,|\kappa|)$
        \State Return $\{a_1,a_2,\dots,a_{\ell}\}$
    \EndProcedure
    
    \Procedure{EquiRippleAlgorithm}{$a=\{a_1,a_2,\dots,a_\ell\}$,$\kappa$}
        \State Set $i=1$ and $a'=\text{None}$
        \State While $a'=\text{None}$ or $a'\neq a$:
        \State \quad $a'=a$
        \State \quad $a=\textsc{Iterate}(a',\kappa)$
        \State \quad $i=i+1$
        \State Return $a$
    \EndProcedure
\end{algorithmic}
\label{alg:iterative}
\end{algorithm}

\begin{figure}[h]
\includegraphics[width=\textwidth]{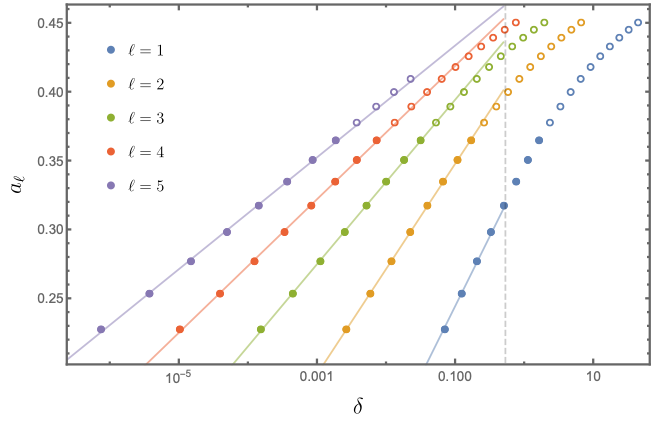}
\caption{Showing the relation between $a_\ell$ and the ripple height $\delta$ for polynomials created by the iterative algorithm. The plot shows best linear fits to the leftmost eight points (filled points) for each $\ell$. The gray vertical line indicates where $\delta=\frac{1}{2}$. Generally, in this case of fixed $\ell$ we observe $\frac{1}{2}-a_\ell\propto\log(1/\delta)$ as $\delta$ gets small. Note that numerical instability begins to set in for polynomials with $\ell=5$ (degree 21) as $a_\ell$ approaches $\frac{1}{2}$ and leads to failure in the $\ell=6$ (degree 25) case for all values of $a_\ell$ we attempted.}
\label{fig:fig1}
\end{figure}

\begin{figure}[h]
    \centering
    \includegraphics[width=\textwidth]{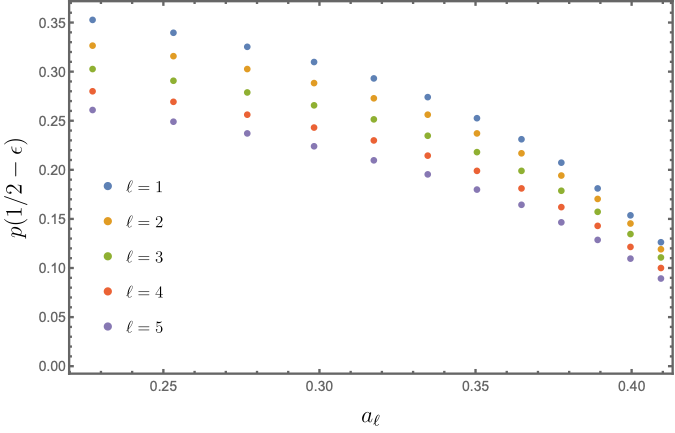}
    \caption{For $\epsilon=\frac{1}{20}$, we show the behavior of $p(\frac{1}{2}-\epsilon)$ as a function of $a_\ell$, where $p$ is the polynomial output from Algorithm~\ref{alg:iterative}. Since empirically $p(\frac{1}{2}-\epsilon)$ is monotonically decreasing with $a_\ell$, and $\delta$ is monotonically increasing with $a_\ell$ (see Figure~\ref{fig:fig1}), we hypothesize there is always a value of $a_\ell$ where $p(\frac{1}{2}-\epsilon)=\delta$.}
    \label{fig:fig2}
\end{figure}


The results of the algorithm are summarized in Figures~\ref{fig:fig1} and \ref{fig:fig2}. Notice that while sometimes the algorithm returns equi-ripple solutions that are not in $\mathcal{P}$ (because the ripple height $\delta$ is larger than 1), one can increase $\ell$ or decrease $a_{\ell}$ such that an equi-ripple solution in $\mathcal{P}$ exists. From the data, we posit the relation $\frac{1}{2}-a_\ell\propto \log(1/\delta)/\ell$. Analyzing the correctness of this algorithm rigorously is a topic of future research.

\bibliographystyle{amsplain}
\bibliography{bibliography}

\end{document}